\documentclass[a4paper]{article}
\pdfoutput=1
\usepackage{graphics}
\usepackage{amsmath,bbm,amssymb,array,verbatim}
\usepackage{hyperref}
\title{Circular Sequences and the Diameter of Multipermutohedra}
\author{Sarang Aravamuthan\thanks{Ignite R\&D Labs, Tata Consultancy Services,  Chennai, India. \newline
\indent \texttt{sarang.aravamuthan@tcs.com}}
}
\newtheorem{thm}{Theorem}

\newtheorem{cor}[thm]{Corollary}
\newtheorem{prop}[thm]{Proposition}

\newcommand{\cp}{\mbox{$\cal P$}}
\newcommand{\cs}{\mbox{$\cal S$}}

\newcommand{\cx}{\mbox{$\cal X$}}

\newcommand{\diam}{\mbox{\rm diam}}
\newcommand{\real}{\hbox{I\kern-.2em\hbox{R}}}

\begin{document}
\maketitle

\begin{abstract}
We derive bounds on the number of switches at an arbitrary set of positions in a 
circular sequence of permutations and relate them to the 
diameter of Multipermutohedra. \newline \newline
\end{abstract}

\section{Introduction}

A sequence of $\bar N +1$ (with $\bar N = {N \choose 2}$)
permutations, $\cp = \langle \pi_0,\ldots,\pi_{\bar N} \rangle $, in the symmetric group $\cs_N$ 
is called a {\it circular sequence} if $\pi_0 = 1\,2\ldots n$ is the identity permutation,
$\pi_{\bar N} = n \, (n-1)\ldots 1$ is the reverse permutation and $\pi_{i+1}$ differs from $\pi_i$ by 
an adjacent transposition, i.e., $\pi_{i+1} = \pi_i\, (j,j+1)$ for some $j$ such that 
$\pi_i(j) < \pi_i(j+1)$. These sequences are used for bounding the number of $k$-sets of a point
configuration \cx~in the plane (a $k$-set of \cx~is a subset of size $k$ separated by
a line from the other points in \cx). The connection between 
$k$-sets and circular sequences was first established by Goodman and Pollack \cite{GP}. For a subset $K \subseteq [(n-1)/2]$, one can also
bound the number of $k$-sets, $k \in K$ over all point configurations \cx. Such bounds were established for specific choices of $K$  using this connection (see \cite{AG, LVWW}).
\smallskip

For the sequence $\cp$, defined above, the process of moving from $\pi_i$ to $\pi_{i+1}$ is referred to as a
{\it switch} $(\pi_i(j), \pi_i(j+1))$ at position $j$ (note that the numbers being swapped are $\pi_i(j)$ and $\pi_i(j+1)$).
We define $s_j(\cp)$ to be the total number of switches at position $j$. 
We count the total number of switches at a given set of positions 
$y = \{y_1,\ldots, y_n\} \subseteq [N-1]$, with $y_1 < y_2 < \cdots < y_n$ which is 
\[ s(\cp,y) = s_{y_1}(\cp) + \cdots + s_{y_n}(\cp). \]

Our goal is to derive a bound on $s(\cp, y)$ and show how it relates to the diameter of \emph{multipermutohedra}
(which, as we will see later, are essentially ``permutohedra on multisets").
\smallskip

In the next section, we estimate $s(\cp, y)$. In Section~\ref{multi}, we introduce the multipermutohedron
and derive bounds on its diameter. In Section~\ref{mpoint}, we show how the 
diameter of multipermutohedra relates to a variant of the $k$-set problem.
\smallskip

We fix some notations. $[n]$ represents the set $\{1,2,\ldots,n \}$.
If $x$ is not an integer, we write $[x]$ instead of $[\lfloor x \rfloor ]$.
For $S \subseteq [n]$ and $x \in \real^n$ define
$x(S) := \sum\limits_{i \in S} x_i$ and denote the size of $S$ by $|S|$.
\smallskip

A \emph{composition} of $n$ is a sequence
of positive integers $\lambda := \langle \lambda_1,\ldots,\lambda_k \rangle $ with $\sum \lambda_i = n$.
\smallskip

Permutations in $\cs_n$ will be represented as words.
We call $e := 1\,2\ldots N \in \cs_N$ the identity permutation and $\bar e := N\,N-1\ldots 1$ the
reverse permutation.
\smallskip

\section{A Bound on Circular Sequences}
\label{circ}

With $\cp$ and $y$ as defined above, we evaluate $s(\cp,y)$ in the following way. 
Each $i \in [N]$ starts at position $i$ in the identity permutation $e$ and
reaches position $N-i+1$ in $\bar e$. Following the notation of \cite{LVWW}, we refer 
to the positions $y_1,\ldots, y_n$ as {\it gates} and
define $c_i$ to be the number of gates that $i$ must cross in order to reach position $N-i+1$.
$c_i$ is interpreted as the cost in moving $i$ to position $N-i+1$.
As $c_i = c_{N-i+1}$, the total cost $c(N)$ is
$$c(N) := \sum_{i=1}^N c_i = 2 \sum_{i=1}^{\lfloor \frac{N}{2} \rfloor} c_i.$$
\smallskip

Next, we define the vector $x = (x_1,\ldots,x_n)$ where $x_j$ represents 
the distance of $y_j$ from one of the ends, i.e.
\begin{equation}
x_j := \min \{y_j, N - y_j\}.
\label{mde1}
\end{equation} 

Let $(p_1,\ldots,p_{n})$ be the permutation of $(x_1,\ldots,x_{n})$ with $N/2 \ge p_1 \ge \cdots \ge p_n > 0$.
Any number $j$ such that $p_{i+1} < j \le p_i$ has cost $c_j = i$ (we assume $p_{n+1} = 0$). Hence
$$ c(N) = 2 \sum_{i=1}^n i(p_i-p_{i+1}) = 2p([n]) = 2x([n]).$$

We now amortize the cost $c(N)$ over $\cp$. Each switch $(i,j)$ across a gate
is interpreted as contributing $+1$ or $-1$ to $c_i$ according to whether $i$ moves towards or
away from position $N-i+1$. A switch $(i,j)$ across a gate is called {\it good} if it contributes $1$ to
both $c_i$ and $c_j$ and {\it bad} if it contributes $1$ to $c_i$ (or $c_j$) and $-1$ to the other.
We observe that every switch across a gate is either good or bad, i.e., must contribute $+1$ to
the cost of at least one of the numbers being moved. This follows from the fact that $i < j$, so if
$i$ is in position $> N - i > N - j$, then $j$ moves towards $N - j + 1$.
\smallskip

Hence, good switches contribute $2$ to $c(N)$ while bad switches contribute $0$.
Since $c(N) = 2 x([n])$, the number of good switches is $x([n])$. Thus if
$s_b(\cp,y)$ is the number of bad switches, then
\begin{equation}
s(\cp,y) = x([n]) + s_b(\cp,y).
\label{mde2}
\end{equation}

Thus, it suffices to estimate $s_b(\cp,y)$, in order to bound $s(\cp,y)$.
Let $l$ be the number of $y_i$ smaller than $N/2$, i.e., $y_l < N/2$
while $y_{l+1} \ge N/2$, and let $r := n - l$. In other words, $l$ and $r$ represent the
number of gates to the $l$eft and $r$ight of $N/2$ respectively. The following result bounds $s(\cp,y)$.

\begin{thm}
Let 
$$s(y) = \min \{s(\cp,y) \, |\:  \cp {\rm \ is\ a\ circular\ sequence} \} $$
be the minimum of the number of switches at a set of positions given by $y$ over all 
circular sequences of permutations in $\cs_N$. Then,
$$ x([n]) \le s(y) \le x([n]) + lr. $$
\label{mdt1}
\end{thm}
\kern-2em

\noindent
{\bf Proof.} 
From (\ref{mde2}), we observe that $s(\cp,y) \ge x([n])$ for any circular sequence $\cp$.
This proves the lower bound for $s(y)$.
\smallskip

To show the upper bound, we construct a circular sequence $\cp$ such that the number of bad switches
$s_b(\cp,y) \le lr$. We construct the sequence in two phases.
\smallskip

In the first phase, we move the numbers $1,\ldots,r$
to positions $y_r,y_{r-1},\ldots y_1$ and the numbers $N,N-1,\ldots,N-l+1$ to positions 
$y_{r+1},\ldots,y_n$ in the following way.
Starting from $e$, we move $1$ to position $y_r$
through a sequence of switches $(1,i), i \le y_r$. 
Next, $2$ is moved to position $y_{r-1}$ with the switches $(2,i), i \le y_{r-1}$.
Continuing this way, the first $r$ and the last $l$ numbers (in the order $N, N-1, \dots, N-l+1$)
are moved to positions $y_r, y_{r-1},\dots, y_1$ and
$y_{r+1}, \dots, y_n$ respectively.

The resulting permutation is
$$ r+1 \ldots r|\ldots r-1|\ldots\ldots 2|\ldots 1|\ldots N|\ldots N-1|\ldots\ldots N-l+1|\ldots N-l$$
where the $|$ indicates the gate positions $y_1,\ldots,y_n$ and the numbers 
$r+1,\ldots,N-l$ are in the remaining positions in increasing order.
\smallskip

The second phase consists of moving $N-i+1$ to position $i$ for $i=1,2,\ldots,\lfloor \frac{N}{2}\rfloor$.
This is done by starting with $N$ and moving it to the first position by switching it in succession
with the numbers $1,\ldots,r$. This brings $1,\ldots,r$ to 
positions $y_{r+1},y_r,\ldots,y_2$ and $r+1$ to position $y_1$ if $y_1 > 1$. Next we move $N-1$ 
to the second position. In general, suppose $\gamma$ is a permutation in this sequence with the 
last $i$ numbers in the first $i$ positions and the first $j$ numbers in the last $j$ positions, 
and $\gamma(i+1) = k$, i.e.
$$\gamma = N\, N-1 \ldots N-i+1\,\, k\, \ldots\, k-1| \ldots k-2| 
\ldots\ldots j+1| \ldots N-i\,\, j\, j-1 \ldots 1.$$
Here the numbers $j+1, \ldots, k-1$ are at the gate positions
and the numbers $k, k+1,\ldots, N-i-1$ are in the
remaining positions in increasing order. Switching $N-i$ in succession with $j+1,j+2,\ldots,k$ moves
$N-i$ to position $i+1$, and $j+1$ to position $N-j$ resulting in the permutation
$$ N\, N-1 \ldots N-i\, \, k+1\ldots k| \ldots k-1| \ldots\ldots j+2| \ldots N-i-1\, \, j+1\, j \ldots 1.$$
Repeating this operation, we finally obtain the permutation $\bar e$.
\smallskip

We observe that all switches in the second phase are good; for
a switch $(i,j)$ across a gate, the numbers $1,\dots,(i-1)$ are already to the right of $i$.
So $i$ is at position $< N-i + 1$. The same argument holds for $j$.
\smallskip
 
In the first phase, suppose $l \le r$. Then the number of bad switches are at most
$ l(r-l)$ for moving the numbers $1,\ldots,r-l$ to positions $y_r,y_{r-1},\ldots,y_{l+1}$; 
$(l-1) + (l-2) + \cdots + 0 = l(l-1)/2 $ for moving the numbers $r-l+1,\ldots,r$ to positions
$y_l,y_{l-1},\ldots,y_1$ and $l + (l-1) + \cdots+1 = l(l+1)/2$ 
for moving the numbers $N, N-1, \ldots, N-l+1$ to positions $y_{r+1},\ldots,y_n$  
adding up to a total of at most $lr$ bad switches. The case $l > r$ is handled similarly. 
This proves the result. 
{\hfill $\Box$}
\medskip

\stepcounter{thm}

\noindent
{\bf Example \thethm}
Let $N = 8$ and $y = \{ 1,4,6, 7\}$. By (\ref{mde1}), $x = (1,4,2,1)$. 
Since $y_1 < N/2$ and $y_2 \ge N/2$, $l = 1$ and $r = n - l = 3$. By Theorem~\ref{mdt1},
$$ 8 = \sum x_i \le s(y) \le 8 + lr = 11$$

To construct a circular sequence $\cp$ with $s(\cp,y) = 11$, we start with the identity permutation
and move the numbers $1,2,8$ to positions $y_3,y_2,y_4$ leading to
$3\, 4\, 5\, 2\, 6\, 1\, 8\, 7$. Next we move $8$ to position $1$ 
by switching it in succession with the numbers $1,2,3$ giving
$8\, 4\, 5\, 3\, 6\, 2\, 1\, 7$. Next $7$ is moved to the second place. 
We summarize by showing some of the permutations in this
sequence.
\begin{eqnarray*}
1\, |2\, 3\, 4\, |5\, 6\, |7\, |8 \quad {\buildrel 4 \over \longrightarrow } & &
\underline{3}\, |4\, 5\, \underline{2}\, |6\, \underline{1}\, |\underline{8}\, |7 
\quad {\buildrel 3 \over \longrightarrow} \\
8\, |\underline{4}\, 5\, \underline{3}\, |6\, \underline{2}\, |\underline{1}\, |\underline{7} 
\quad {\buildrel 3 \over \longrightarrow} & &8\, |7\, \underline{5\, 4}\, |\underline{6}\, 3\, |2\, |1 
\quad {\buildrel 1 \over \longrightarrow}  \quad 8\, 7\, 6\, 5\, 4\, 3\, 2\, 1.
\end{eqnarray*}
where the numbers above the arrows indicate the number of switches required to move to the next 
permutation, the `$|$' show the positions where the switches are counted and 
the numbers being switched are underlined. {\hfill $\Box$}

\newcounter{temp}
\setcounter{temp}{\thethm}
\medskip

We see that the lower bound for $s(y)$ is attained when all the gates are to the left (or right) of the middle. In this
case, $r = 0$ (resp. $l = 0$) and the two bounds for $s(y)$ coincide. The upper bound is attained when $y = [N-1]$.
In this case, $s(y) = {N \choose 2}$, i.e. we count switches at all positions.

\section{Circular Sequences and Multipermutohedra}
\label{multi}

The multipermutohedron is a generalization of the permutohedron $P_N$ which is the 
convex hull of all permutations of the point $(1,2,\ldots,N) \in \real^N$.
A multipermutohedron is defined by taking the convex hull of all permutations of a {\it multiset}, 
that is, a ``set'' with repeated elements. 
\smallskip

Let $b_1 < b_2 < \cdots < b_n $ be $n$ distinct numbers ($n > 1$) and consider the multiset with
$k_i$ copies of $b_i$ for $i = 1, \dots, n$, i.e.
\begin{equation}
M := \{ b_1^{k_1},b_2^{k_2},\ldots,b_n^{k_n}\} = \{a_1,\ldots,a_N\} \hskip 1cm  k_i > 0,\  i=1,\ldots,n  
\label{me2}
\end{equation}
with $ a_1 \le a_2 \le \cdots \le a_N$ and $N = \sum\limits_{i=1}^n k_i$.
Let $P(M)\subset \real^N$ be the polytope formed by
taking the convex hull of
all permutations of the point $(a_1,\ldots,a_N) \in \real^N$. 
We call $P(M)$ the {\it multipermutohedron} on the multiset $M$.
\smallskip

Multipermutohedra were first studied by Schoute \cite{Sc}, who constructed them all from the 
simplex by a sequence of simple operations of {\it expansion} and {\it contraction}. 
Independent treatments of multipermutohedra also appear in \cite{BS2, Sch}. The
proofs for the assertions that follow can be found in \cite{BS2}.
\smallskip

The inequality description of $P(M)$ is a direct generalization of the
description for the permutohedron.
$$P(M)=\{x\in \real^N \ | \ x([N])  =  \sum_{i=1}^N a_i, \
x(S)  \ge  \sum_{i=1}^{|S|} a_i,\ \hbox{\rm  for all}\ S\subset [N]\ \}. $$
Thus if $M = \{1,2,2,3\}$, then $P(M) \subset \real^4$ is given by
\begin{eqnarray*}
x_1 + x_2 + x_3 + x_4 & = & 8 \\
1 \le x_i & \le & 3 \hskip 1cm i = 1,\ldots,4 \\
x_i + x_j & \ge & 3 \hskip 1cm 1 \le i < j \le 4.
\end{eqnarray*}

\begin{figure}[htbp]
\begin{center}
\includegraphics{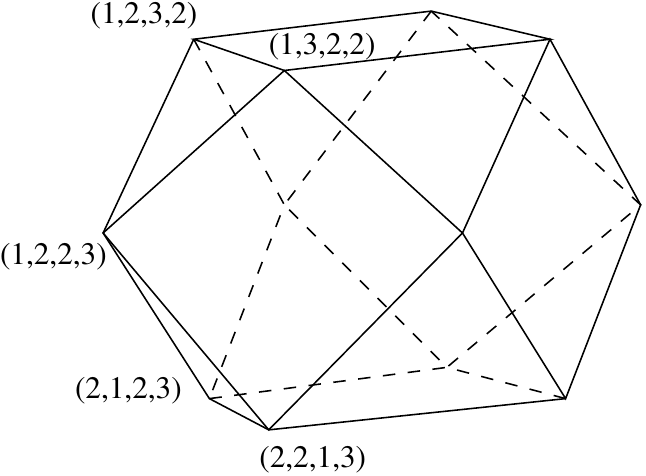}
\end{center}
\caption{A $3$-dimensional multipermutohedron $P(\{1,2,2,3\})$.}
\label{mfig1}
\end{figure}

Figure~\ref{mfig1} shows $P(\{1,2,2,3\})$.  It has 14 facets corresponding to
these 14 inequalities, which correspond, in turn, to the 14 ordered
partitions of the set $\{1,2,3,4\}$ into two parts. 
\smallskip

Faces of $P(M)$ are products of lower dimensional multipermutohedra and
correspond to ordered partitions of $[N]$. As in the permutohedron, the face lattice of the 
multipermutohedron does not depend on the numbers being permuted but only on their multiplicities.
Henceforth we assume that, in (\ref{me2}), $b_i =i$ and rewrite $M$ as 
\begin{equation}
M = \{1^{k_1},\ldots,n^{k_n} \}.
\label{MMM}
\end{equation}
\smallskip

Given $\sigma \in \cs_N$, the vertex $a_\sigma$ has $1$'s in positions 
$\sigma(1),\ldots,\sigma(k_1)$, $2$'s in positions 
$\sigma(k_1+1),\ldots, \sigma(k_1+k_2)$ and so on, i.e.
\begin{equation} 
a_\sigma := (a_{\sigma^{-1}(1)},\ldots,a_{\sigma^{-1}(N)}).
\label{me3}
\end{equation}

Note that this correspondence is not one-to-one; several permutations will correspond 
to the same vertex of $P(M)$.
\smallskip

Adjacency of vertices on $P(M)$ is obtained by switching two components whose values are consecutive.
The relation to permutations is as follows. For $j=1,\ldots,n$, define
\begin{equation}
y_j := k_1+\cdots +k_j
\label{me4}
\end{equation}
and let $y_0 = 0$.
We define {\it blocks} of consecutive integers $S_1,\ldots,S_n$ by
$$ S_j := [y_j]\setminus [y_{j-1}].$$
Given $\sigma,\pi \in \cs_N$,
the vertices $a_\sigma$ and $a_\pi$  
are adjacent if $a_\pi = a_{\sigma (i,j)}$ for some transposition $(i,j)$ with $i$ and $j$ in
successive blocks. In other words, $\pi$ is obtained from $\sigma$ by switching $\sigma(i)$ and $\sigma(j)$ 
for some $i$ and $j$ in successive blocks and possibly permuting the values of $\sigma (i,j)$ 
within each block.
\medskip

\stepcounter{thm}

\noindent
{\bf Example \thethm }
Let $M = \{1,2,2,3\}$ as earlier. Let $e = 1\, 2\, 3\, 4$ be the identity permutation.
This corresponds to the vertex $a_e = (1,2,2,3)$ of $P(M)$. Multiplying $e$ by transpositions $(i,j)$
with $i,j$ in successive blocks leads to the permutations 
$2\, 1\, 3\, 4$, $3\, 2\, 1\, 4$, $1\, 4\, 3\, 2$ and $1\, 2\, 4\, 3$.
These in turn correspond to the vertices $(2,1,2,3)$, $(2,2,1,3)$, $(1,3,2,2)$ 
and $(1,2,3,2)$ that are adjacent to $(1,2,2,3)$ as shown in Figure \ref{mfig1}.
{\hfill $\Box$}
\medskip

For $\sigma,\pi \in \cs_N$, we denote by $d(\sigma,\pi)$, the shortest distance between the 
vertices $a_\sigma$
and $a_\pi$ in $P(M)$. For convenience, we also denote the vertex $a_\sigma$ in $P(M)$ by $\sigma$. 
\smallskip

The \emph{diameter} of $P(M)$ (denoted diam($P(M)$)) is the largest of the distances between 
pairs of vertices of $P(M)$.

\begin{prop}
The vertex farthest from $\sigma \in P(M)$ is 
${\bar \sigma} = \sigma(N)\,\sigma(N-1) \ldots \sigma(1)$. 
The diameter of $P(M)$ is $d(e,\bar e)$.
\label{mdp0}
\end{prop}
{\bf Proof.} 
We first show that the vertex farthest from $e$ in $P(M)$ is ${\bar e}$.
This will follow if we can show that for any permutation $\pi$ with 
$a_\pi \not= a_{\bar e}$, switching a pair of numbers $\pi(i) < \pi(j)$ with $i< j$
in successive blocks yields a permutation $\pi' = \pi (i,j)$ with $d(e,\pi') \ge d(e,\pi)$.
\smallskip

Suppose $d(e,\pi) > d(e,\pi')$. Let $i' = \pi(i)$ and $j'=\pi(j)$. Since $\pi$ and $\pi'$ are 
adjacent vertices in $P(M)$, $d(e,\pi) = d(e,\pi')+1$. Consider the shortest path
from $e$ to $\pi$ through $\pi'$. Let $\gamma, \alpha_1$ be successive vertices on this path with
$\alpha_1^{-1}(j')< \alpha_1^{-1}(i')$ in successive blocks while $\gamma^{-1}(i')$ and $\gamma^{-1}(j')$
are either in the same block or in successive blocks. 
Hence, in going from $\gamma$ to $\alpha_1$, we have moved $j'$ to the left of $i'$.
Let the vertices on the path from $\gamma$ be
$\gamma,\alpha_1,\ldots,\alpha_{p-1} = \pi', \alpha_p = \pi$. We construct a shorter
path from $\gamma$ to $\pi$, namely $\gamma, \beta_1,\ldots,\beta_{p-1} = \pi$ 
where for each $m$, $\beta_m = (i',j')\alpha_m$. We see that for each $m > 0$, $\beta_{m+1}$ 
differs from $\beta_m$
by a transposition across adjacent blocks showing that they are adjacent. 
Since $\beta_1$ switches $i'$ and $j'$ in $\alpha_1$,
we see that either
$\beta_1 = \gamma$ or $\beta_1$ and $\gamma$ are adjacent (this case subsumes the possibility of $\beta_1 = \alpha_1$).
In any case we see that 
this path is shorter than our chosen path by at least $1$.
This contradicts our hypothesis that $d(e,\pi) > d(e,\pi')$ proving our claim.
\smallskip

Given a path from $e$ to $\pi$ in $P(M)$, multiplying the permutations 
in this path by $\sigma$ yields a path from $\sigma$ to $\sigma \pi$ of the same length. 
Hence the vertex farthest from $\sigma$ is $\sigma \bar e = \bar \sigma$ and 
$d(\sigma,\bar \sigma) = d(e,\bar e)$ proving the result.
{\hfill $\Box$}
\medskip

We now relate the diameter of $P(M)$ to circular sequences. Let the multiset $M$ and 
$y = (y_1,\ldots,y_{n-1})$ be given by (\ref{MMM}) and (\ref{me4}) respectively.
From the criterion for adjacency of vertices on $P(M)$, we observe that each circular sequence
$\cp$ corresponds to a path 
in $P(M)$ from $e$ to $\bar e$ of length $s(\cp,y)$. From Proposition~\ref{mdp0}, it follows that
$$\diam\, (P(M)) = \min \{s(\cp,y) \, |\:  \cp {\rm \ is\ a\ circular\ sequence} \}. $$

\noindent
Theorem~\ref{mdt1} automatically translates to bounds on diameter of $P(M)$.

\begin{cor}
Let $M$ be defined by (\ref{MMM}) with $k_i > 0$ and let $y = (y_1,\ldots,y_{n-1})$ and 
$x = (x_1,\ldots,x_{n-1})$ be specified by (\ref{me4}) 
and (\ref{mde1}). 
Then the diameter of $P(M)$ is bounded by 
$$ x([n-1]) \le \diam\, (P(M)) \le x([n-1]) + lr $$
with $l$ and $r$ as defined in Theorem~\ref{mdt1}. {\hfill $\Box$}
\label{multi2}
\end{cor}

\stepcounter{thm}

\noindent
{\bf Example \thethm } Consider the multiset
$M = \{1, 2^3, 3^2,4,5 \}$. The vectors $y$ and $x$ given by (\ref{me4}) and (\ref{mde1}) are
$y = (1,4,6,7)$ and $x = (1,4,2,1)$. Since $y_1 < N/2$ and $y_2 \ge N/2$, $l=1$ and $r = n-1-l = 3$.
By Corollary~\ref{multi2},
$$ 8 = \sum x_i \le \diam\, (P(M)) \le 8 + lr = 11.$$

A circular sequence $\cp$ with $s(\cp,y) = 11$ was constructed in Example~\arabic{temp}
and this translates to a path of length $11$ between $e$ and $\bar e$ in $P(M)$. {\hfill $\Box$}
\medskip

It's easy to see that the upper bound in Corollary~\ref{multi2} is attained for the permutohedron $P_N$ which
has a diameter of ${N \choose 2}$. Also, for small $n$, it's possible to derive an explicit expression for diam($P(M)$).

\begin{prop}
Let the multiset $M$ be given by the composition $\langle k_1,\ldots,k_n\rangle $ of $N$. If $n=2$ then
$\diam\, (P(M)) = \min \{k_1,k_2 \}$. When $n=3$,
$$ \diam\, (P(M)) = 
\begin{cases} x_1 + x_2 & \mbox{\rm if } k_1 \not= k_3;\\  x_1+x_2+1 = 2k_1+1 & \mbox{\rm if } k_1 = k_3.
\end{cases}
$$
\label{mdp2}
\end{prop}
\noindent
{\bf Proof.}  
For $n=2$, the lower and upper bounds for the diameter in Corollary~\ref{multi2} 
are the same and the result follows.
If $n=3$ and $k_1 \not= k_3$ then it is easy to show by a careful choice of switches that 
the lower bound for the diameter is attainable. If $k_1 = k_3$, then by Corollary~\ref{multi2}, 
the diameter is either $2k_1$ or $2k_1+1$. The lower bound is not attained because, for any 
circular sequence, the switches at positions $k_1$ and $N-k_1$
cannot all involve a number less than $k_1$ and a number greater than $N-k_1$, i.e., there must be
at least one bad switch. This proves the result. {\hfill $\Box$}
\medskip

We observe that since the diameter of $P(M)$ is obtained by counting 
the number of switches at certain positions
in a circular  sequence in $\cs_N$, its value is at most the total number of switches i.e. $N \choose 2$. 
It would be interesting to determine, in some form, all the integers in the set $[{N \choose 2}]$
that are the diameters of multipermutohedra $P(M) \subset \real^N$ for some $M$ given by (\ref{MMM}).

\section{Multipermutohedra and \texorpdfstring{$k$}{k}-sets}
\label{mpoint}

We now relate the diameter of the multipermutohedron to arrangements of points on a plane.
Let $\cx \subset \real^2$ be a configuration of $N$ points in general position on a plane,
 i.e., no three points of $\cx$ lie on a line.
For $k \le \lfloor N/2 \rfloor$, we define a {\it left} 
(resp. {\it right}) {\it $k$-set} to be a set of $k$ points of $\cx$ that
lie on the left (resp. right) of a line (with respect to a directed reference line, say the X-axis). 
For a line parallel to the X-axis, we take the left of the line to be the open half-space above the line.
A $k$-set of $\cx$ is either a left $k$-set or a right $k$-set. Let $f_l(k,\cx)$, $f_r(k,\cx)$ and
$f(k,\cx)$ denote the number of left $k$-sets, right $k$-sets and $k$-sets of $\cx$ respectively. 
\smallskip

For subsets $L,R \subseteq [(N-1)/2]$, we count the number of 
sets that appear as a left $k$-set for $k \in L$ or a right $k$-set for $k \in R$. We define

\begin{equation}
f(L,R,\cx) := \sum_{k \in L \cap R} f(k, \cx) + \sum_{k \in L\setminus R} f_l(k, \cx) + 
\sum_{k \in R \setminus L} f_r(k, \cx).	       
\label{mdpe1}
\end{equation}

When $L = R \subseteq [(N-1)/ 2]$, we are counting the number of $k$-sets for $k \in L$
and we write $f(K,\cx)$ for $f(K,K,\cx)$. Our objective is to derive bounds for $f(L,R,\cx)$ in terms of the
diameters of certain multipermutohedra. 

Since the minimum (resp. maximum) of $f(L,R,\cx)$ is not affected by a slight perturbation 
of the points of $\cx$, we assume that no two points of $\cx$ lie on a line parallel to the $X$-axis.
Then, for each $k \le (N-1)/2$, there are two $k$-sets
that are both left $k$- and right $k$-sets. Hence
\begin{equation}
f_l(k, \cx) + f_r(k,\cx) = f(k,\cx) + 2. 
\label{mdpe0}
\end{equation}

Following the approach of Goodman and Pollack \cite{GP}, we associate with $\cx$
a circular sequence of permutations in the following way.
Project the points of $\cx$ onto a directed line $\ell$ parallel to the Y-axis. 
Label the projected points on $\ell$ in the order $1,\ldots,N$ and rotate $\ell$ counter clockwise.
The projection of the points on $\ell$ determines a permutation $\pi \in \cs_N$.
The first $k$ points $\pi(1),\ldots, \pi(k)$ form a left $k$-set 
while the last $k$ points $\pi(N-k+1),\ldots, \pi(N)$ form a right $k$-set. When $\ell$ becomes
perpendicular to a line joining a pair of points of \cx, the order of the projection 
of these two points on $\ell$ is
reversed. The new permutation differs from $\pi$ by an adjacent transposition.
When $\ell$ has rotated through an angle of $180^o$, the order of the projected points becomes
$N,N-1,\ldots,1$. Thus $\cx$ describes a circular sequence of permutations 
$\cp(\cx) = \langle e,\pi_1,\ldots,\pi_{\bar N} = \bar e \rangle $.
If $s_k$ is the number of switches at position $k$ in $\cp(\cx)$,
then the number of left $k$-sets $f_l(k, \cx) = s_k+1$ and the number of 
right $k$-sets $f_r(k,\cx) = s_{N-k}+1$. Hence by (\ref{mdpe0}), the number of $k$-sets 
$f(k,\cx) = s_k + s_{N-k}$. By (\ref{mdpe1}),
\begin{equation}
f(L,R,\cx) = \sum_{k \in L} s_k + \sum_{k \in R} s_{N-k} + |L \cup R| - |L \cap R|.
\label{mdpe2}
\end{equation}
The subsets $L$ and $R$ define a composition $\langle L,R \rangle = \langle k_1,\ldots,k_n\rangle $ of $N$ 
where $n = |L| + |R|+1$
and the elements of $L$ are the partial sums $k_1 + k_2 + \cdots + k_i$ that are at most $(N-1)/2$ 
while the elements of
$R$ are the partial sums $k_n + k_{n-1}+ \cdots + k_j$ that are at most $(N-1)/2$. 
Let $P(\langle  L,R \rangle )$ be the multipermutohedron
defined by this composition. Then the sequence $\cp(\cx)$ describes a 
path in $P(\langle  L,R \rangle )$ from $e$ to $\bar e$ of
length $\sum\limits_{k \in L} s_k + \sum\limits_{k \in R} s_{N-k}$ which must be at least its diameter.
This gives a lower bound for $f(L,R,\cx)$.
\smallskip

To bound $f(L,R,\cx)$ from above, we consider the subsets $\bar L := [N/2]\setminus L$ and 
$\bar R := [(N-1)/2]\setminus R$.
As before, these define a composition
$\langle \bar L, \bar R \rangle$ of $N$ and a multipermutohedron $P(\langle  \bar L, \bar R \rangle )$. 
Since $\sum\limits_{k = 1}^{N-1} s_k = {N \choose 2}$, we can rewrite (\ref{mdpe2}) as 
$$f(L,R,\cx) = {N \choose 2} -\left( \sum_{k \in \bar L} s_k +
\sum_{k \in \bar R} s_{N-k} \right) +|L\cup R| -|L \cap R|.$$
This gives the upper bound for $f(L,R,\cx)$ and ties the diameter of $P(M)$ to $k$-sets 
in the following way.

\begin{thm}
Let $\cx$ be a configuration of $N$ points in $\real^2$ in general position 
and let $L,R \subseteq [(N-1)/2]$. Then the function $f(L,R,\cx) - |L \cup R| + |L \cap R|$ is bounded from 
above and below by the diameter of multipermutohedra
$P(\langle L,R\rangle )$ and $P(\langle \bar L,\bar R\rangle )$, that is,
$$\diam\, (P(\langle  L,R \rangle ))  \le f(L,R,\cx) - |L \cup R| + |L \cap R| \le 
{N \choose 2} - \diam\, (P(\langle  \bar L, \bar R \rangle )) . $$
In particular,
$$ \diam\, (P(\langle  K,K \rangle )) \le f(K,\cx) \le 
{N \choose 2} - \diam\, (P(\langle  \bar K,\bar K \rangle ))$$
thus bounding the number of $k$-sets for $k \in K \subseteq [(N-1)/2]$ by diameters of 
multipermutohedra.
{\hfill $\Box$} 
\label{mdpt1}
\end{thm}

The problem of $k$-sets and its relation to circular sequences has been extensively studied 
(see \cite{AG,ELSS,GP,LVWW}). In \cite{LVWW}, $k$-sets are used to derive a lower bound on the number
of convex quadrilaterals in a set of $n$ points in the plane.
\smallskip

A lower bound for the number
of $k$-sets of a point configuration in general position follows from Proposition~\ref{mdp2} 
and the above theorem. As observed in \cite[Example 8]{LVWW}, this lower bound of $2k+1$
can actually be achieved by a configuration of $N$ points ($N \ge 2k+1$) which consist
of a regular ($2k+1$)-gon with the remaining points situated close to the center of the gon. 
Finding the upper bound for the number of $k$-sets
seems to be a harder problem. Some estimates for this bound are given in \cite{ELSS,TD}.
\medskip

\noindent
{\bf Acknowledgment.} The author would like to thank L.~J.~Billera for helpful discussions.

\end{document}